\documentclass[12pt,a4paper]{amsart}
\usepackage{amssymb}
\usepackage{latexsym}
\usepackage{amsthm}

\begin{document}

\theoremstyle{plain}

\newtheorem{thm}{Theorem}[section]
\newtheorem{cor}[thm]{Corollary}
\newtheorem{lem}[thm]{Lemma}
\newtheorem{prop}[thm]{Proposition}
\newtheorem{ex}[thm]{Example}
\newtheorem{ax}[thm]{Axiom}
\newtheorem{ques}[thm]{Question}
\newtheorem{defn}{Definition}[section]
\newtheorem{rem}{Remark}[section]

\title{Arithmetic of the [19,1,1,1,1,1] fibration}
\author{Matthias Sch\"{u}tt}
\address{Institut f\"{u}r Mathematik (C), Universit\"{a}t
Hannover, Welfengarten~1, 30167 Hannover, Germany}
\email{schuett@math.uni-hannover.de}
\author{Jaap Top}
\address{IWI-RuG, P.O.Box 800, 9700 AV Groningen, 
the Netherlands}
\email{top@math.rug.nl}
\date{\today}
\thanks{It is a pleasure to thank Jasper Scholten for his
interest in this work and his verification of part of the table
below. The first author gratefully acknowledges partial support 
from the 
DFG-Schwerpunkt 1094 ``Globale Methoden in der komplexen 
Geometrie''.}
\maketitle
\section{Introduction}
Suppose $\pi:S\rightarrow {\mathbb P}^1$ is an elliptic $K3$-surface over ${\mathbb C}$
with a section, such that all its singular fibres are
multiplicative. Suppose moreover that $\pi$ is extremal,
i.e., the group of sections of $\pi$ is finite and the
rank $\rho$ of the N\'eron-Severi group $NS(S)$ of $S$ is maximal,
so $\rho=20$. A formula of Shioda and Tate \cite[Cor.~5.3]{ShMW}
then implies that $\pi$ has precisely $6$ singular fibers,
consisting of $n_i$ components ($i=1,\ldots,6$) and $\sum n_i=24$.
A well-known result of Miranda and Persson \cite{MP} gives
all $112$ realizable $6$-tuples $[n_1,n_2,n_3,n_4,n_5,n_6]$.

The largest number $n_i$ appearing in this table is
$n_i=19$. In fact, T.~Shioda \cite{ShCR} showed for
 an
elliptic $K3$-surface over ${\mathbb C}$ with a section that the property of 
containing
a multiplicative fibre consisting of $19$ components,
determines it uniquely. Such a fibration was found
explicitly by Y.~Iron (a student of R.~Livn\'{e}) in his
master's thesis \cite{Ir} and also by Shioda \cite{ShCR}, who 
realized the
relation with work of M.~Hall,~jr. \cite{Hall}.
The result can be stated as follows.
\begin{prop}[Iron, Shioda]\label{one}
Consider the polynomials 
\[f=t^8+6t^7+21t^6+50t^5+86t^4+114t^3+109t^2+74t+28
\] and
\[\begin{array}{r}
g=t^{12}+9t^{11}+45t^{10}+156t^9+408t^8+846t^7+1416t^6
+1932t^5\\+2136t^4+1873t^3+(2517t^2+1167t+299)/2.\end{array}
\]
Let $S$ be the elliptic $K3$-surface corresponding to the
Weierstrass equation 
\[
y^2=x^3-27f(t)x+54g(t).
\]

Then $S\to{\mathbb P}^1$ has precisely $6$ singular fibres.
The fibres over the zeroes of $4t^5+15t^4+38t^3+61t^2+62t+59$
are of type $I_1$. At $t=\infty$, the elliptic curve over
${\mathbb Q}(t)$ defined by the given Weierstrass equation has
split multiplicative reduction. The fibre of
$S\to{\mathbb P}^1$ over $t=\infty$ is of type $I_{19}$.
\end{prop}

The fact that the reduction at $t=\infty$ is split
multiplicative, is easily verified. Indeed,
the change of variables $\xi:=x/t^4, \eta:=y/t^6, s:=1/t$
yields a Weierstrass equation 
where the special fibre at $s=0$,
is given as $\eta^2=\xi^3-27\xi+54=(\xi-3)^2(\xi-3+9)$.
This clearly defines a fibre with split multiplicative
reduction. It follows that all $19$ components of the
fibre in $S$ over $t=\infty$ are defined over ${\mathbb Q}$. 
Since by the formula of 
Shioda and Tate \cite[Cor.~5.3]{ShMW}, these components together 
with the zero-section 
generate the N\'{e}ron-Severi group of $S$, we obtain
\begin{cor}
The $K3$-surface $S$ over ${\mathbb Q}$ defined in
Proposition~\ref{one} satisfies
$\rho(S/{\mathbb Q})=20$, where $\rho(S/{\mathbb Q})$ denotes
the rank of the subgroup $NS(S/{\mathbb Q})$ of the
N\'{e}ron-Severi group $NS(S)$ generated by the classes
of ${\mathbb Q}$-rational divisors.
\end{cor}

It should be noted that this corollary contradicts
\cite[Thm.~1]{ShRanks}. Indeed, the argument presented in
{\em loc.cit.} asserts that for almost all prime numbers $p$,
the reduction $S\bmod p$ satisfies $\rho(S/{\mathbb F}_p)
=20$ or $22$. We will show that this is incorrect. More 
precisely, on the one hand, $\rho(X/{\mathbb F}_p)=22$ is 
impossible for $p\neq 2$ and $X/{\mathbb F}_p$ a $K3$-surface,
 due to an argument of M.~Artin \cite[(6.8)]{Artin}. On the 
other hand, we will see in the following (Cor.~\ref{cor}) that 
$\rho(S/{\mathbb F}_p)=21$ for a set of primes with 
density 1/2. 
As a consequence, the proof of Corollary~3 in \cite{ShRanks}
is incomplete. This means that the following question is
still open:
\begin{ques}
Is the maximal rank $rk\,E({\mathbb Q}(t))$ for an
elliptic curve $E/{\mathbb Q}(t)$ corresponding to an
elliptic $K3$-surface over ${\mathbb Q}$ equal to $18$
or smaller than $18$?
\end{ques}

Recently, Shioda published an erratum \cite{Sh05} to his paper
\cite{ShRanks}. It turns out that several experts had pointed
out the given mistake, including K.~Hulek and H.~Verrill.
Interestingly, Shioda's original method of proof 
shows that $\rho(A/{\mathbb Q})\neq 4=h^{1,1}(A_{\mathbb C})$,
 when $A/{\mathbb Q}$ is an abelian surface. 
In the rest of the present note, we will describe 
(Section~\ref{2}) the
representation of the Galois group 
$\mbox{Gal}(\overline{\mathbb Q}/{\mathbb Q})$
on $H^2(S_{\overline{{\mathbb Q}}},{\mathbb Q}_{\ell})$.
This allows us to compute $\rho(S/{\mathbb F}_p)$
and the Hasse-Weil zeta function of $S/{\mathbb F}_p$,
for almost all primes $p$, supplementing the remarks
made by Shioda in \cite[\S~2]{Sh05}. 
In Section~\ref{3}, we recall a conjecture which Shioda
formulated in \cite[\S~4]{ShCR}.
We use the surface $S$ to illustrate how this conjecture can be
verified in specific cases.

We remark as an aside, that similar examples may be found
starting from other surfaces in the Miranda-Persson list
\cite{MP}. Explicit equations for some of these
surfaces can be obtained from \cite{TY} and \cite{Sch},
and more recently for all of these surfaces from \cite{BM}.

\newpage

\section{$H^2$ as a Galois module}\label{2}
Throughout, $S$ is the elliptic surface over ${\mathbb Q}$
defined in Proposition~\ref{one}. Let $G_{\mathbb Q}:=
\mbox{Gal}(\overline{\mathbb Q}/{\mathbb Q})$.
As a $G_{\mathbb Q}$-module, the ${\mathbb Q}_{\ell}$
vector space $H^2(S_{\overline{{\mathbb Q}}},{\mathbb Q}_{\ell})$
splits as a direct sum into the two-dimensional representation 
of the transcendental lattice and one-dimensional 
representations for the N\'eron-Severi group:
\[
H^2(S_{\overline{{\mathbb Q}}},{\mathbb Q}_{\ell})=
T_{\ell} \oplus \oplus_{j=1}^{20} {\mathbb Q}_{\ell}(1).
\]
Here ${\mathbb Q}_{\ell}(1)$ denotes the one dimensional
${\mathbb Q}_{\ell}$ vector space on which a Frobenius
conjugacy class at $p\neq \ell$ acts as multiplication
by $p$. Our aim is to obtain the trace and the determinant
of Frobenius conjugacy classes on $T_{\ell}$.

For this purpose we have to restrict to the primes of good 
reduction, which in our case
include all primes $p\neq 2, 3, 19$. We now discuss the
three remaining primes.

For $p=3$, a model for $S$ over ${\mathbb Q}$ 
with good reduction can be derived as follows. 
Consider the ${\mathbb Q}$-isomorphic surface coming
from the equation 
\[y^2=x^3-\frac{1}{3}f(t) x+\frac{2}{27} g(t)\]
and apply the change of variables 
$x\mapsto x+(t^4+3t^3+3t^2+4t+4)/3$. This gives rise to the 
equation
$y^2=x^3+A(t)x^2-B(t)x+C(t)$ with 
\[\begin{array}{l}
A=t^4+3t^3+3t^2+4t+4,\\
B=2t^6+8t^5+14t^4+22t^3+23t^2+14t+4\\
C=t^8+5t^7+12t^6+21t^5+27t^4+24t^3+15t^2+5t+1.
\end{array}\]
The corresponding surface is easily seen to have good reduction 
at 3. 

Concerning the reduction at 19, 
Shioda \cite[Thm. 3.1]{ShCR} shows that the minimal 
desingularization of the nodal Weierstrass surface over 
$\overline{\mathbb{F}}_{19}$ requires one further blow-up 
compared to characteristic 0. 
Reasoning differently, we will see below that the
representation of $G_{\mathbb Q}$ on $T_\ell$ ramifies
at $19$. Hence the reduction at 19 is bad.

At $p=2$, it turns out that no elliptic fibration on a
$K3$-surface in characteristic two containing an $I_{19}$-fibre
exists. This follows from a direct calculation, by adopting
techniques of W.E.~Lang \cite[\S~4]{La} to the present
situation. The details may be found in \cite{Sch2}.

We shall now compute the determinant of a Frobenius conjugacy class on $T_{\ell}$ for a good prime $p\neq 2, 19$.
From \cite[\S~3]{ShCR} (or by a direct calculation using
the intersection matrix on $NS(S)$), one derives that
the intersection form on the transcendental lattice of $S$
is given by ${2\;\;1}\choose{1\;\;10}$. This implies
(compare \cite[Example~1.6]{L}) that the determinant
of $T_{\ell}$ is given by 
\[
\det(\mbox{Frob}_p)=p^2\epsilon(p).
\]
Here $\epsilon$ is the quadratic character corresponding
to the extension\linebreak ${\mathbb Q}(\sqrt{-19})/{\mathbb Q}$.
Explicitly, $\epsilon(p)=1$ when
$p\equiv 1, -2, -3, 4, 5, 6, 7, -8$, or $9 \bmod 19$ and 
$\epsilon(p)=-1$ for all other
primes $p\neq 19$.

The trace of Frobenius elements acting on $T_{\ell}$ is
related to the number of points in $S({\mathbb F}_p)$ 
via the Lefschetz trace formula. Using the description of $H^2$ given above,
this says in our situation that
\[
\#S({\mathbb F}_{p^n})=
1+20p^n+
\mbox{trace}\left((\mbox{Frob}_p)^n|T_{\ell}\right)+p^{2n}.
\]
Write $q=p^n$. Of the ${\mathbb F}_q$-rational points in $S$,
precisely $19(q+1)-19=19q$ are on the $19$ components
over $t=\infty$. The remaining ones are on the fibres
over $t\in{\mathbb F}_q$, hence correspond to
rational points on the projective cubic given by
$y^2z=x^3-27f(t)xz^2+54g(t)z^3$. (For $p=3$, we have to work with the other model given above.) This gives $q$ rational points
for $z=0$. All other points are obtained from the affine
Weierstrass equation $y^2=x^3-27f(t)x+54g(t)$. Their
number is
\begin{small}
\[
\sum_{x,t\in{\mathbb F}_q} \left(1+\chi_q(x^3-27f(t)x+54g(t))\right)
=q^2+\sum_{x,t\in{\mathbb F}_q}\chi_q(x^3-27f(t)x+54g(t))\]\end{small}  
where $\chi_q:{\mathbb F}_q\to\{0,\pm 1\}$ sends $0$ to $0$
and is the unique quadratic character on ${\mathbb F}_q^*$ for $p\neq 2$.
Combining the two formulas above yields
\[
\mbox{trace}\left((\mbox{Frob}_p)^n|T_{\ell}\right)=-1+
\sum_{x,t\in{\mathbb F}_q}\chi_q(x^3-27f(t)x+54g(t)).
\]
Here is a small table computed using this formula. For convenience, we also include the trace of $(\mbox{Frob}_p)^2$, although we will not need it in the following.

\vspace{\baselineskip}

\begin{center}
\begin{tabular}{|r|r|r|}\hline 
$p$ & trace $\mbox{Frob}_p$ on $T_{\ell}$ &
 trace $(\mbox{Frob}_p)^2$ on $T_{\ell}$\\\hline
 3 & 0 & -18\\
  5 & -9 & 31\\
 7 & -5 & -73\\
  11 & 3 & -233\\
13 & 0 & -338\\
17& 15 & -353\\
23 & -30 & -158\\
29& 0 &-1682\\
31&0& -1922\\
37& 0 & -2738\\
41 & 0 & -3362\\
43 & -85 & 3527 \\\hline
\end{tabular}
\end{center}

\vspace{\baselineskip}

According to \cite[1.3 -- 1.5]{L}, the $L$-series
determined by $T_{\ell}$ equals the $L$-series of
a cusp form of weight $3$ and nebentypus the quadratic
character $\epsilon$ of conductor $19$. We will now relate this cusp form to a Hecke character.

Let $K={\mathbb Q}(\sqrt{-19})$ and denote by ${\mathbb A}_K^*$
the ideles of $K$.
Then
\[
\psi_0(x):=x_{\infty}^{-2}\prod_{v\;\mbox{\scriptsize finite}}
\pi_v^{2v(x_v)}
\]
(in which $\pi_v$ generates the maximal ideal in $O_v\cap K$)
defines the unique unramified Hecke character $\psi_0:{\mathbb A}_K^*/K^*\to
{\mathbb C}^*$ with $\infty$-type 2.
Using the table above, an argument completely analogous (in fact
somewhat simpler since the field $K$ has class number one) to
\cite[\S~2]{L} shows that the representations of
$G_{\mathbb Q}$ given by $\mbox{Ind}_{G_K}^{G_{\mathbb Q}}\psi_0$
and by $T_{\ell}$ agree. 
Alternatively, this could be checked directly via the
Faltings-Serre method \cite{Se}; in the present case, it
suffices to check that the characteristic polynomials of
Frobenius at the primes $3, 5, 7$ and $31$ agree.

Explicitly, this implies the following.
\begin{prop}
Let $p\neq 2, 19$ be a prime number. The characteristic
polynomial of $\mbox{Frob}_p$ acting on $T_{\ell}$ is
given by one of these rules:
\begin{enumerate}\label{two-one}
\item If $p$ is inert in ${\mathbb Q}(\sqrt{-19})$, then this
polynomial equals $X^2-p^2$.
\item If $p$ splits in ${\mathbb Q}(\sqrt{-19})$, then write
$4p=a^2+19b^2$ for integers $a,b$. Put $c:=(a^2-19b^2)/2$.
Then the required characteristic polynomial is
$X^2-cX+p^2$.
\end{enumerate}
\end{prop} 
The description of the characteristic polynomials here follows
from the definition of the Hecke $L$-series $L(s,\psi_0)$.
Its Euler factor at an inert prime of $K$ is
$\left(1-\pi_v^2p^{-2s}\right)^{-1}=(1-p^{2-2s})^{-1}$.
The corresponding characteristic polynomial is $X^2-p^2$.
In case $p$ splits in $K$, one can write $p=\pi\bar{\pi}$
for some $\pi=(a+b\sqrt{-19})/2$ in the ring of integers
of $K$. The Euler factor now is
$\left((1-\pi^2p^{-s})(1-\bar{\pi}^2p^{-s})\right)^{-1}=
(1-cp^{-s}+p^{2-2s})^{-1}$. This yields the other
characteristic polynomials. $\qed$

\begin{cor} \label{cor}
Let $p\neq 2, 19$ be a prime number. 
One has 
\[\rho(S/{\mathbb F}_p)=\left\{\begin{array}{cl}
20& \mbox{if  } p
 \mbox{  splits in  } {\mathbb Q}(\sqrt{-19});
\\ 21 & \mbox{if  } p>2 \mbox{  is inert in  } {\mathbb Q}(\sqrt{-19}).
\end{array}\right.
\]
\end{cor}
Indeed, for $p>3$ this follows from the Tate conjecture. Since it is known for elliptic 
$K3$-surfaces over
finite fields of characteristic $p>3$ \cite{Ta}, the rank
$\rho(S/{\mathbb F}_p)$ equals the multiplicity of the
factor $X-p$ in the characteristic polynomial of $\mbox{Frob}_p$
acting on $H^2(S_{\overline{{\mathbb F}}_p},{\mathbb Q}_{\ell})$. For the verification in case $p=3$, we refer to the following example.
$\qed$

\vspace{\baselineskip}

\begin{ex}\label{2.3}
Consider the reduction $S/{\mathbb F}_3$ at 3 as constructed
 above. Our aim is to determine $NS(S/\overline{\mathbb F}_3)$. 
Since the singular fibre configuration in $S\to {\mathbb P}^1$
is the same in characteristic $3$ as in characteristic $0$, 
the Tate conjecture predicts $S/\overline{\mathbb F}_3$ to have 
Mordell-Weil group of rank 2 with the sections defined over 
${\mathbb F}_3$ as a rank $1$ subgroup.   
We  find a section over ${\mathbb F}_3$ as 
\[P=(t^2+t, t-1).\] 
A further computation gives another section 
\[Q=(-t^2+(1+i)t+1, t^4+it^3+(1+i)t-i)\]
over ${\mathbb F}_9={\mathbb F}_3[i]$ where $i^2=-1$. In order 
to prove that these generate the Mordell-Weil group
of the given fibration $S\to{\mathbb P}^1$ over 
$\overline{\mathbb F}_3$, we use the height pairing as defined 
by Shioda in \cite{ShMW}. An easy calculation 
(compare Section~\ref{3}) gives
\[
<P,P>=\frac{6}{19},\;\;\;\;\;<P,Q>=
\frac{9}{19},\;\;\;\;\;<Q,Q>=\frac{42}{19}.
\] 
As a consequence, the sublattice $N$ of 
$NS(S/\overline{\mathbb F}_3)$ which is generated by $P, Q$ and 
the trivial lattice consisting of the zero-section and the 
components of the fibres has discriminant 
\[
 (-19)\det\left(\frac{3}{19}\begin{pmatrix}
2 & 3\\ 3 & 14
                  \end{pmatrix}\right)=-9.
\]
But by the general theory \cite[p.~556]{Artin}, a supersingular 
$K3$-surface has 
discriminant $-p^{2\sigma_0}$ with $\sigma_0\in\{1,\hdots,10\}$ 
the Artin invariant. Hence $N=NS(S/\overline{\mathbb F}_3)$ as 
claimed. In particular, the Tate conjecture holds for $S$ in
characteristic $3$.
\end{ex}

It is well known (see, e.g., \cite[Section~2.4]{To89}) that a
$2$-dimensional Galois representation as described here comes
from a modular form. In the present case, this is the
unique normalized newform of weight $3$ for $\Gamma_0(19)$ with
nebentypus the quadratic character $\epsilon$ of conductor $19$,
which has Fourier coefficients in ${\mathbb Z}$. 
Here, the uniqueness is shown as follows.
Since the weight is odd, a newform with totally real
Fourier coefficients has CM by its own nebentypus 
\cite[Proposition 3.3]{Ribet}.
If the level of the modular form equals 
the discriminant of the CM field (up to sign), then 
the corresponding Hecke character has trivial conductor.
If moreover the class number of the CM field is one,
this Hecke character and
hence the modular form is uniquely determined.
Alternatively, in the present case the uniqueness may be
read off from the tables \cite{St}.

Since all eigenvalues of Frobenius on $H^2$ are explicitly
given, one obtains the local zeta function 
$Z(S/{\mathbb F}_p,T)$.
\begin{prop}
\begin{enumerate}
\item For primes $p>2$ which are inert in ${\mathbb Q}(\sqrt{-19})$
one has
\[
Z(S/{\mathbb F}_p,T)=\frac{\displaystyle 1}{\displaystyle
(1-T)(1+pT)(1-pT)^{21}(1-p^2T)}.
\]
\item For primes $p$ which split in ${\mathbb Q}(\sqrt{-19})$
one has
\[
Z(S/{\mathbb F}_p,T)=\frac{\displaystyle 1}{\displaystyle
(1-T)(1-cT+p^2T^2)(1-pT)^{20}(1-p^2T)},
\]
with $c$ as defined in Proposition~\ref{two-one}(2).
\end{enumerate}
\end{prop}

\section{Reductions modulo supersingular primes}\label{3}
Modulo every odd prime $p$ which is non-split in 
${\mathbb Q}(\sqrt{-19})$, the surface $S$ defines a
supersingular elliptic $K3$-surface over ${\mathbb F}_p$.
For these primes, there is an injective reduction map
$NS(S/{\mathbb Q})\to NS(S/\overline{{\mathbb F}}_p)$
from a rank $20$ lattice to a lattice of rank $22$.
The orthogonal complement of the image will be denoted
by $L(p)$; it is a negative definite even lattice of
rank two.
In \cite[\S~4]{ShCR}, Shioda formulates the conjecture
that such a lattice $L(p)$ should be similar to the
transcendental lattice $T_S$ of $S/{\mathbb Q}$.
Recall that the intersection form on $T_S$ is
given by ${2\;\;1}\choose{1\;\;10}$.
Here we will verify this conjecture for the reductions
of $S$ modulo $3$ and modulo $19$.

In characteristic $3$, two extra generators of
$NS(S/\overline{\mathbb F}_3)$ are given by the points
$P$ and $Q$ defined in Example~\ref{2.3}.
Using a calculation as in the proof of \cite[Lemma~8.1]{ShMW}
 one finds that 
(in Shioda's notations)
\[\varphi(P):=(P)-(O)-2F+\sum_{i=1}^{4}
\frac{14i}{19}\Theta_i+
\sum_{i=5}^{18} \frac{5(19-i)}{19}\Theta_i
\]
and
\[
\varphi(Q):=(Q)-(O)-2F+\frac{2}{19}\Theta_1+
\sum_{i=2}^{18} \frac{17(19-i)}{19}\Theta_i
\]
generate the vector space $L(3)\otimes {\mathbb Q}$.
Note that from this, the height pairings as
given in Example~\ref{2.3} can be checked.
In this vector space, $L(3)$ is generated by $19\varphi(Q)$
and $\varphi(P)-7\varphi(Q)$.
As a result, the intersection form on $L(3)$ can be given by
$-3 {{\;266\;\;-95}\choose{-95\;\;\;34}}$.
This lattice is easily checked to be similar to $T_S$.

In characteristic $19$, the equation $y^2=x^3-27f(t)x+54g(t)$
becomes 
\[
y^2=x^3-8(t+3)(t-5)^7x-3(t+3)^{11}(t-5).
\]
It follows from \cite[Theorem~3.1]{ShCR} that the Mordell-Weil
group is torsion-free of rank one, and is generated by
\[
R:=\left(2\frac{(t+3)^{10}}{(t-5)^6},
7\sqrt{-1}\frac{(t+3)^{15}}{(t-5)^9}\right).
\]
Apart from the $I_{19}$-fibre over $t=\infty$, the
elliptic fibration has one other reducible fibre:
over $t=-3$ there is a fibre of type $III$.
Its component not meeting the zero-section we denote by $\Theta$.
The lattice $L(19)$ is generated by $\Theta$ and $(R)-(O)-5F$.
On this basis, the intersection form is
$-{{\;2\;\;\;-1}\choose{-1\;\;10}}$.
Hence also in this case, Shioda's conjecture holds.

Note that, in the notations of \cite{ShMW},
$\varphi(R)=(R)-(O)-5F+\frac12\Theta$.
From this it can be seen, that the lattice generated
by $\Theta$ and $2\varphi(R)$ is 
$\langle -2\rangle\oplus \langle -38\rangle$,
which is {\em not} similar to
the transcendental lattice.

\end{document}